\documentclass[3p,twocolumn]{elsarticle} 
\usepackage{lineno,hyperref,xcolor}
\usepackage{rotating}
\usepackage{setspace}
\usepackage{amsthm,amsmath,amsfonts}
\hypersetup{colorlinks=true,linkcolor=blue,filecolor=magenta,urlcolor=blue,citecolor=blue,}
\biboptions{authoryear}

\begin{document}
\begin{frontmatter}
\title{A stopping rule for randomly sampling bipartite networks with fixed degree sequences}
\journal{Social Networks}
\author{Zachary P. Neal}
\ead{zpneal@msu.edu}
\address{Michigan State University}

\begin{abstract}
Statistical analysis of bipartite networks frequently requires randomly sampling from the set of all bipartite networks with the same degree sequence as an observed network. Trade algorithms offer an efficient way to generate samples of bipartite networks by incrementally `trading' the positions of some of their edges. However, it is difficult to know how many such trades are required to ensure that the sample is random. I propose a stopping rule that focuses on the distance between sampled networks and the observed network, and stops performing trades when this distribution stabilizes. Analyses demonstrate that, for over 650 different degree sequences, using this stopping rule ensures a random sample with a high probability, and that it is practical for use in empirical applications.\\
\\
{\textcolor{red}{\textbf{This is a post-print. The version of record is available at: Neal, Z. P. (2024). A stopping rule for randomly sampling bipartite networks with fixed degree sequences. \textit{Social Networks, XX}, XXX. \url{https://doi.org/10.1016/j.socnet.XXX}}}}
\end{abstract}

\end{frontmatter}


\section{Introduction}
Statistical analysis of bipartite networks frequently requires randomly sampling bipartite networks with fixed degree sequences \citep[e.g.,][]{neal2021comparing,jasny2012baseline,gotelli2000null}. The sampled networks provide a comparison (i.e. a null model) against which to compare and thereby evaluate the statistical significance some network property of interest, such as the nestedness of a bipartite network \citep{strona2014methods} or the strength of an edge in a bipartite projection \citep{neal2014backbone}. Trade algorithms offer an efficient way to iteratively randomize a bipartite network while preserving its degree sequences by trading the positions of edges \citep{strona2014fast,godard2022fastball}. If enough trades are performed, these algorithms are known to yield random samples of all such networks \citep{carstens2015proof}. However, the number of trades required to ensure that the sample is random remains unknown.

In this paper, I propose a stopping rule for these algorithms that focuses on the distribution of sampled networks' distances from the observed network. Using this rule, trades are performed until this distribution has stabilized. Examining this stopping rule's performance in over 650 degree sequences, I demonstrate that on average it yields a random sample over 85\% of the time.

The paper is organized in four sections. In the background section I review trade algorithms for randomizing and sampling bipartite networks, illustrating their performance in a small bipartite network. I then propose a stopping rule for deciding how many trade randomizations to perform. In the methods section, I describe a test of its performance for bipartite networks with varying degree sequences. In the results section, I present the results of this test, and illustrate the stopping rule's practicality using empirical examples drawn from ecology, sociology, and political science. Finally, in the discussion section, I summarize the findings, and discuss the method's limitations and implications.

\section{Background}
\subsection{Sampling bipartite networks using trade randomization}
A network is \textit{bipartite} when its nodes can be partitioned into two sets, such that edges exist only between nodes in different sets. A bipartite network is also \textit{two-mode} when these two sets of nodes represent distinctly different kinds of entities. Such networks are common in the literature, and can capture many kinds of phenomena including event attendance where people are connected to events \citep{breiger1974duality}, bill co-sponsorship where legislators are connected to bills \citep{neal2014backbone}, ecosystems where species are connected to habitats \citep{gotelli2000null}, and scholarship where authors are connected to papers \citep{newman2001scientific}. Following \cite{latapy2008basic}, I generically refer to the smaller set of nodes as `top nodes,' and to the larger set of nodes as `bottom nodes.'  A bipartite network is characterized by two degree sequences, one for the top nodes, and another for the bottom nodes. For example, in a bipartite network capturing event attendance, where people are top-nodes and events are bottom-nodes, the top-node degree sequence captures the number of events attended by each person, while the bottom-node degree sequence captures the number of people attending each event.

\begin{figure}
    \centering
    \includegraphics[width=\linewidth]{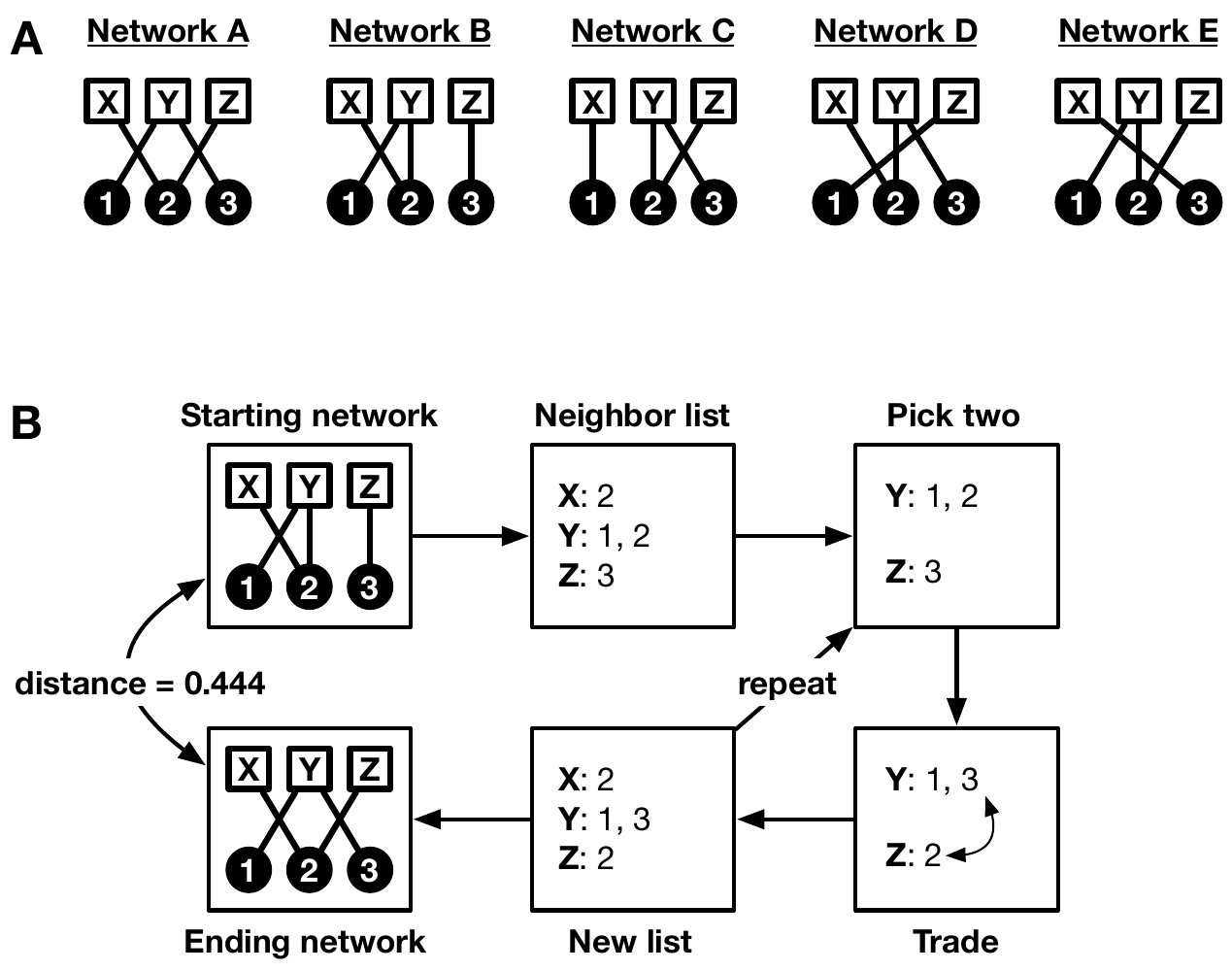}
    \caption{(A) The set $\mathcal{B}$ of all bipartite networks with top-degree sequence \{1,2,1\} and bottom-degree sequence \{1,2,1\}; (B) Example of bipartite network randomization using a trade algorithm.}
    \label{fig:algorithm}
\end{figure}

For any pair of degree sequences, there exists a (possibly empty) set $\mathcal{B}$ of all bipartite networks that have the given degree sequence. For example, given a top node degree sequence \{1,2,1\} and a bottom node degree sequence \{1,2,1\}, $\mathcal{B}$ contains the five networks shown in Figure \ref{fig:algorithm}A, and therefore the size or cardinality of $\mathcal{B}$ is five (i.e., $|\mathcal{B}| = 5$). The goal of randomly sampling bipartite networks with fixed degree sequences is to sample networks from $\mathcal{B}$ such that each member of $\mathcal{B}$ has an equal probability of being sampled. The degree sequences \{1,2,1\} and \{1,2,1\} define a set that is small enough that it is possible to enumerate all its members, making it easy to directly sample from this list. However, because longer degree sequences define sets that are too large to enumerate, and indeed whose sizes are generally unknown \citep{barvinok2010number}, random sampling in such cases requires a different approach.

Many different approaches for sampling have been proposed \citep{verhelst2008efficient,miller2013exact,admiraal2008networksis,bezakova2007sampling,chen2006sequential,fout2020non,wang2020fast,strona2014fast,godard2022fastball}, however they vary in the evidence of their randomness, their computational complexity, and their practicality in empirical contexts. The most widely used sampling approach involves starting with a network (typically, the observed network being analyzed) that has the given degree sequences, then repeatedly rewiring its edges in a way that preserves the degree sequence, to obtain a new network that can be regarded as randomly sampled from $\mathcal{B}$. Among the most efficient algorithms proven to perform such rewiring randomly are the `trade' algorithms \citep{strona2014fast,godard2022fastball,carstens2015proof}, so-called because they trade edges in the same way that children might trade baseball cards. 

Figure \ref{fig:algorithm}B illustrates each step of the trade algorithm. First, the \textit{starting network} (e.g., Network B from Figure \ref{fig:algorithm}A) is represented as a \textit{neighbor list} that indicates, for each top-node, the bottom-nodes to which it is connected. Next, the algorithm \textit{picks two} top-nodes at random from the neighbor list (e.g., nodes Y and Z). Then, the two chosen top-nodes \textit{trade} a random number of their unique neighbors. A single `trade' can involve swapping none, one, or many neighbors, which makes trade algorithms more efficient than checkerboard swap algorithms that swap only one neighbor in each iteration \citep{ellison2024switching}. The outcome of the trade is represented as a \textit{new list}. The process can be repeated, performing additional trades among pairs of randomly-selected top-nodes. Two variants of the trade algorithm exist that differ only in how efficiently they perform these steps: The `curveball' \citep{strona2014fast} algorithm performs trades in \textit{O(n log n)} time, while the newer `fastball' \citep{godard2022fastball} algorithm performs them slightly more efficiently in \textit{O(n)} time. After performing the chosen number of trades, the resulting neighbor list represents an \textit{ending network}. In this example, performing one trade transformed the starting network into a new network with the same degree sequences (e.g., Network A from Figure \ref{fig:algorithm}A), or put another way, resulted in sampling Network A from the set of five possible networks.

\begin{figure*}
    \centering
    \includegraphics[width=\linewidth]{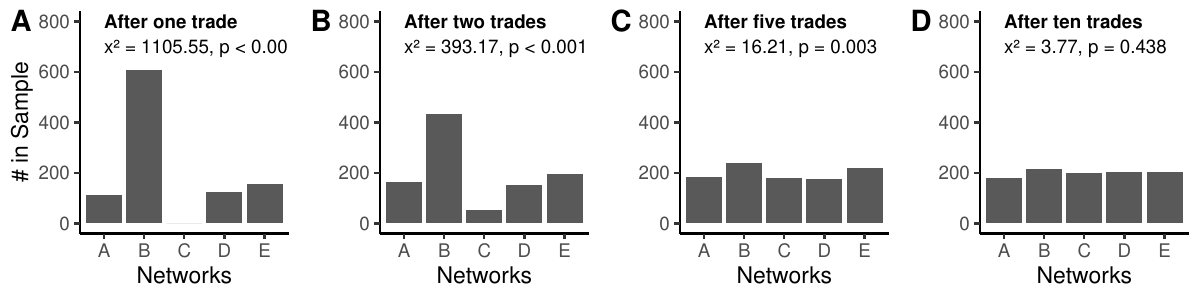}
    \caption{Distribution of networks in a sample drawn using a trade algorithm, following the given number of trades, with a $\chi^2$ test of uniform distribution.}
    \label{fig:composition}
\end{figure*}

Figure \ref{fig:composition} illustrates the composition of a sample of 1000 networks drawn using the fastball trade algorithm, starting from Network B. Performing only one trade is insufficient to obtain a random sample from $\mathcal{B}$ because, as Figure \ref{fig:composition}A shows, it yields a sample in which Network B is over-represented and Network C is under-represented. The statistically significant $\chi^2$ statistic confirms that there is significant variation in the relative frequency of each network in the sample. This variation remains even after two (Figure \ref{fig:composition}B) or five (Figure \ref{fig:composition}C) trades have been performed. However, after ten trades have been performed (Figure \ref{fig:composition}D), each network occurs roughly the same number of times. The uniformity of the sample is confirmed by a non-significant one-sample $\chi^2$ statistic, which indicates that the sample is now a random sample from $\mathcal{B}$. 

Prior work on trade algorithms has focused on proving that they yield random samples from $\mathcal{B}$ if a sufficient number of trades are performed \citep{carstens2015proof}, or on developing ways to perform those trades more quickly \citep[e.g.,][]{godard2022fastball,carstens2018speeding}. However, no prior work has directly investigated how many such trades are required to ensure random sampling (i.e., the mixing time), which is unknown and ``intractable to compute'' \citep[][p. 5]{carstens2015proof}. When \cite{strona2014fast} initially developed the curveball trade algorithm, they provided implementations of their algorithm in R and Python, in which they fixed the number of trades performed at $5T$, where $T$ is the number of top nodes. However, they did not provide a rationale for choosing this value, or evidence supporting its adequacy for yielding random samples. Therefore, in the next section, I propose a stopping rule to decide when a sufficient number of trades have been performed.

\subsection{A trade randomization stopping rule}
As Figure \ref{fig:composition} illustrates, when $|\mathcal{B}|$ is small, a one-sample $\chi^2$ test can be used to determine when enough trades have been performed to obtain a random sample from $\mathcal{B}$. However, typically $|\mathcal{B}|$ is too large to enumerate all its members, and the $\chi^2$ test cannot be used to verify a sample's randomness or to know when to stop performing trades. Therefore, a different method is required for deciding when a sufficient number of trades have been performed to obtain a random sample.

One alternative focuses not on the frequency of each network in the sample, but on the distance between each sampled network and the starting network. The distance $d$ between two bipartite networks is defined as the fraction of dyads that are different, or when the networks are represented as binary incidence matrices, the fraction of cells that are different \citep{strona2014fast}. Formally 
$$d_{A,B} = \frac{\sum{|A-B|}}{r \times c},$$
where $A$ and $B$ are binary $r \times c$ incidence matrices. For example, Figure \ref{fig:algorithm} illustrates that the distance between the starting network (Network B) and the ending network (Network A) is $d_{A,B} = 0.444$. 

\begin{figure*}
    \centering
    \includegraphics[width=\linewidth]{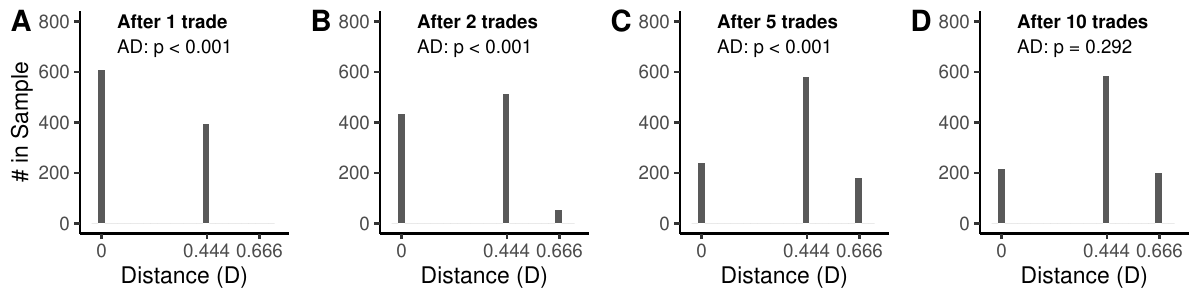}
    \caption{Distribution of distance between a starting matrix and networks in a sample ($D$) drawn using a trade algorithm, following the given number of trades, with an Anderson-Darling (AD) test of equality with an earlier distribution.}
    \label{fig:distance}
\end{figure*}

Given a sample of bipartite networks obtained by randomizing a starting network using a trade algorithm, $D_t$ is the list of distances $d$ between the starting network and each network in the sample after performing $t$ trades. Figure \ref{fig:distance} illustrates the distribution of $D$ in each sample whose composition is shown in Figure \ref{fig:composition}. For example, Figure \ref{fig:distance}A shows that, after one trade, the sample contains many networks at a distance of 0 from the starting network (these are Network B), some networks at a distance of 0.444 (these are Networks A, D, and E), and no networks at a distance of 0.666 (these would have been Network C). Note that this pattern mirrors the frequency of networks in the sample shown in Figure \ref{fig:composition}A: many instances of Network B, some instances of Networks A, D, and E, and no instances of Network C.

As more trades are performed, the distribution of $D$ changes. Specifically, because trade algorithms randomly sample from $\mathcal{B}$ \citep{carstens2015proof}, the distribution of $D$ converges on the distribution of distances $\mathcal{D}$ between the starting network and \textit{every} member of $\mathcal{B}$. This convergence is visible in Figure \ref{fig:distance}, and in this small example can be verified because $\mathcal{D}$ is known: \{0,0.444,0.444,0.444,0.666\}. When the distribution of $D$ matches the distribution of $\mathcal{D}$, this suggests that the members of the sample are representative of, and thus a random sample from, $\mathcal{B}$.

Typically $\mathcal{D}$ is unknown, so it is not possible to directly test whether the distribution of $D$ matches the distribution of $\mathcal{D}$. However, because $D$ converges on $\mathcal{D}$, it is only necessary to test when the distribution of $D$ has stabilized.  There are many ways to test whether $D_t$ and $D_{t-1}$ were drawn from the same distribution. The Kolmogorov–Smirnov test \citep{kolmogorov1933sulla,smirnov1948table} is the most widely known such test. However, because this test lacks sensitivity in a range of cases, two other tests are more appropriate in this context. When the values in $D$ are unequally spaced (i.e., ordinal; as they are in the example illustrated in Figure \ref{fig:distance}), the Anderson-Darling (AD) test \citep{anderson1952asymptotic} is optimal \citep{dowd2020}. In contrast, when the values in $D$ are equally spaced (i.e., interval), the DTS test is optimal \citep{dowd2020}.

The AD test reported in Figure \ref{fig:distance}B shows that the distribution of $D$ after two trades is statistically significantly different from the distribution of $D$ after one trade; the distribution has not stabilized. In contrast, the AD test reported in Figure \ref{fig:distance}D shows that the distribution of $D$ after ten trades is not statistically significantly different from the distribution of $D$ after five trades; after ten trades the distribution has stabilized. The stability of the distribution of $D_{10}$ suggests that the sample is random after 10 trades, which was confirmed directly by the $\chi^2$ test reported in Figure \ref{fig:composition}D.

Therefore, \textit{I hypothesize that sampling bipartite networks with fixed degree sequences by performing trades until the distribution of $D_t$ is not statistically significantly different from the distribution of $D_{t-n}$ will ensure with high probability that the resulting sample is a random sample from $\mathcal{B}$}.

\section{Method}
\subsection{Validation}
The small example above illustrates one set of degree sequences (i.e., \{1,2,1\}, \{1,2,1\}) where, when the distribution of $D_t$ became not statistically significantly different from the distribution of $D_{t-n}$ (see Figure \ref{fig:distance}D), the sample of networks was a random sample from $\mathcal{B}$ (see Figure \ref{fig:composition}D). I follow a similar approach to test the hypothesis that the proposed stopping rule yields random samples, evaluating how often it yields a random sample when used to sample bipartite networks with other degree sequences. Specifically, I consider 665 different degree sequences that define $\mathcal{B}$ with $|\mathcal{B}| \leq 1000$.

For each degree sequence, I apply the stopping rule to draw a sample of 1000 bipartite networks, then evaluate whether the resulting sample is a random sample from $\mathcal{B}$. I repeat this 100 times, recording the percentage of samples that are random. Applying the stopping rule requires computing $D$ and performing a AD or DTS test after every $n$ trades. To avoid needing to perform these computations too often, $n$ is equal to the number of top nodes. The statistical significance of the test is evaluated using the conventional $p < 0.05$ threshold.

To determine whether a given sample is random sample from $\mathcal{B}$, following the approach illustrated in Figure \ref{fig:composition}, I use a one-sample $\chi^2$ test. A sample is uniformly random when this test is not statistically significant, indicating that the distribution of matrix frequencies does not differ from a uniform distribution.

\subsection{Application}
The validation test evaluates whether the stopping rule consistently ensures a random sample, but it is only possible for shorter degree sequences that define $\mathcal{B}$ with smaller $|\mathcal{B}|$. To examine the stopping rule's practicality for randomly sampling from larger spaces, I illustrate its use to draw samples of bipartite networks in three different contexts: ecology, sociology, and political science. In each case, I draw a sample of 1000 random bipartite networks because this is the typical size of a sample required to evaluate the significance of patterns in the network of interest \citep{strona2014fast}.

First, I consider the `Darwin's Finches' network, which records whether Darwin observed each of 13 species of finches living on each of 17 islands \citep{sanderson2000testing}. It is defined by the top-node degree sequence \{14,13,14,10,12,2,10,1,10,11,6,2,17\}, which describes the number of islands on which each species of bird was observed, and by the bottom-node degree sequence \{4,4,11,10,10,8,9,10,8,9,3,10,4,7,9,3,3\}, which describes the number of species that were observed on each island. In the ecological context, randomly sampling bipartite networks with these degree sequences is necessary to identify, for example, species' significant co-habitation on islands that suggests their symbiosis or lack of competition. Interestingly, \cite{chen2005sequential} computed that $|\mathcal{B}| = 67,149,106,137,567,626$ for these degree sequences. Thus, even this modestly-sized bipartite network highlights the large size of $\mathcal{B}$ from which networks are to be randomly sampled, and thus the difficulty in determining whether the sample is random.

Second, I consider the `Southern Women' network, which records whether each of 18 women attended each of 14 social events \citep{homans1950,breiger1974duality}. It is defined by the top-node degree sequence \{8,7,8,7,4,4,4,3,4,4,4,6,7,8,5,2,2,2\}, which describes the number of events each woman attended, and by the bottom-node degree sequence \{3,3,6,4,8,8,10,14,12,5,4,6,3,3\} who attended each event. In the sociological context, randomly sampling bipartite networks with these degree sequences is necessary to identify, for example, women's significant co-attendance at events that suggests their friendship.

Finally, I consider the `US Senate' network, which records whether each of 100 US Senators sponsored each of 5357 in the 117\textsuperscript{th} session of the U.S. Senate \citep{fowler2006legislative,neal2022constructing}. It is defined by a top-node degree sequence that describes the number of bills each Senator sponsored, and by a bottom-node degree sequence that describes the number of Senators who sponsored each bill. In the political science context, randomly sampling bipartite networks with these degree sequences is necessary to identify, for example, Senators' significant co-sponsorship of bills that suggests their political alliance.

Each of these empirical examples involve a $\mathcal{B}$ that is too large to use a $\chi^2$ test to determine whether a sample obtained using the stopping rule is random. Instead, in these illustrative applications, I focus on three outcomes that are observable. First, does the distribution of $D$ converge to a stable distribution that is expected to match the distribution of $\mathcal{D}$? Second, is each member of the resulting sample unique, suggesting that the sample provides broad coverage of $\mathcal{B}$ and does not over-represent certain members of $\mathcal{B}$. Third, can a likely-to-be-random sample be drawn in a practical amount of time using the stopping rule?

\section{Results}
\subsection{Validation}\label{sec:small}
Table \ref{tab:results} reports, for 10 selected degree sequences and for the average over all 665 degree sequences, the percent of samples that were random when the stopping rule is used to decide when to stop performing trade randomizations. Each row represents a pair of bipartite network degree sequences, which define a $\mathcal{B}$ containing a specific number of networks ($|\mathcal{B}|$). For example, the first row corresponds to the small example used throughout the background section. In this small example, using the proposed stopping rule to decide when to stop performing trade randomizations yielded a random sample of bipartite matrices from $\mathcal{B}$ 92\% of the time. Using the stopping rule for other degree sequences yielded a random sample with similar frequency, and on average over 365 different degree sequences it yielded a random sample 87.5\% of the time. 

\begin{table*}[]
\caption{Performance of stopping rule for selected degree sequences}
\label{tab:results}
\begin{tabular}{lllll}
\hline
\multicolumn{2}{c}{\underline{Degree Sequences}} & & Mean Trades & Random \\
Top Nodes & Bottom Nodes & $|\mathcal{B}|$ & Required & Samples \\
\hline
\{1,2,2\} & \{1,2,2\} & 5 & 10 & 92\% \\
\{2,2,3\} & \{1,2,2,2\} & 12 & 10 & 95\% \\
\{1,3,3\} & \{1,1,1,2,2\} & 18 & 10 & 88\% \\
\{1,4,4\} & \{1,1,1,2,2,2\} & 24 & 11 & 91\% \\
\{1,2,5\} & \{1,1,1,1,1,1,2\} & 51 & 11 & 90\% \\
\{1,3,6\} & \{1,1,1,1,1,1,2,2\} & 72 & 11 & 93\% \\
\{2,2,2,3\} & \{1,2,3,3\} & 27 & 15 & 92\% \\
\{2,2,2,2\} & \{1,1,1,2,3\} & 108 & 14 & 96\% \\
\{1,1,3,4\} & \{1,1,1,1,2,3\} & 54 & 18 & 89\% \\
\{1,2,3,4,4\} & \{1,3,3,3,4\} & 76 & 24 & 84\% \\
... & ... & ... & ... & ... \\
\multicolumn{2}{l}{\textbf{Average for 665 sequences}} & \textbf{172} & \textbf{17} & \textbf{87.50\%}\\
\hline
\multicolumn{5}{l}{Results for all 665 degree sequences are provided in \textit{Supplementary Information}.}
\end{tabular}
\end{table*}

These estimated rates of achieving random samples must be interpreted with caution because they are likely underestimated. The one-sample $\chi^2$ test used to determine whether a sample is uniformly random is performed with $\alpha = 0.05$, which gives the test a type-I error rate of 5\%. This means that the test will incorrectly reject the null hypothesis (i.e., find that the sample is not random) despite the null hypothesis being correct (i.e., the sample actually is random) in 5\% of cases. Therefore, even if the stopping rule yielded a random sample 100\% of the time, the $\chi^2$ test would find it did so only 95\% of the time. Adjusting these rates to account for this type-I error suggests that the stopping rule's average success rate is 92\% (i.e., $\frac{0.875}{0.95}$).

Whether the raw or adjusted success rate is used, these findings provide supports the hypothesis: the proposed stopping rule ensures with high probability that the resulting sample is a random sample from $\mathcal{B}$.

\subsection{Application}\label{sec:large}
\begin{figure*}
    \centering
    \includegraphics[width=\linewidth]{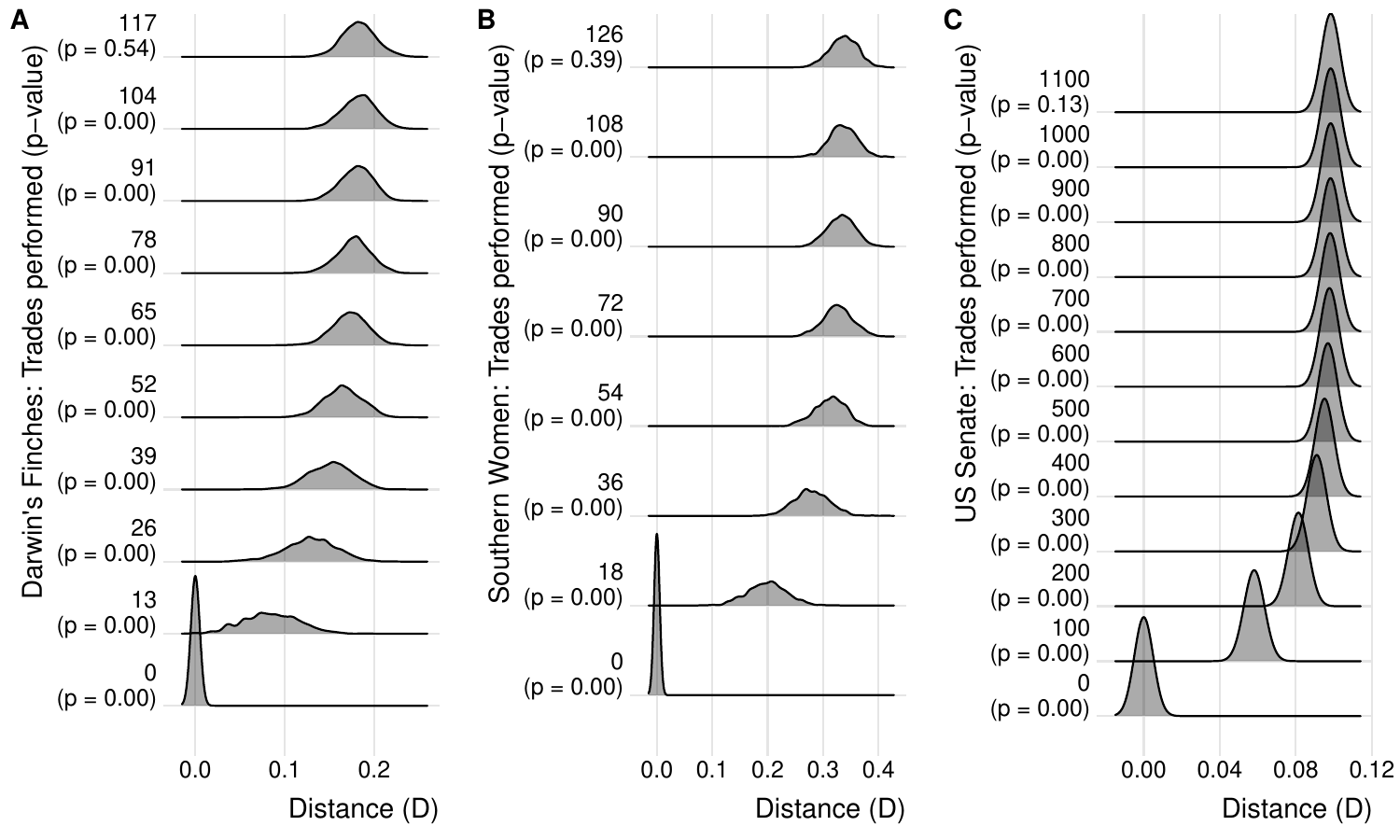}
    \caption{Distribution of sampled networks' distance from observed (A) Darwin's Finches, (B) Southern Women, and (C) US Senate networks following trades}
    \label{fig:example}
\end{figure*}

Figure \ref{fig:example} shows the distribution of $D$, after a specified number of trade randomizations, in a sample of bipartite networks with the same degree sequences as the (A) Darwin's Finches, (B) Southern Women, and (C) US Senate networks. In each case, the stopping rule stopped performing trade randomizations when this distribution stabilized, which required a different number of trades depending on the network (e.g., 117 trades in the finches network, 1100 trades in the Senate network). It is impossible to verify that the resulting samples are random samples from the corresponding $\mathcal{B}$. However, it is possible to examine whether the stopping rule performs as expected and in practical time.

First, as expected, in all three cases the distribution of $D$ converges to a stable distribution. Notably, the distribution on which it converges is unique to the degree sequence, which is consistent with the expectation that the distribution of $D$ converges on the distribution of $\mathcal{D}$.

Second, as required of a random sample, all 1,000 members of these samples are unique. This suggests that applying the stopping rule yields a sample with broad coverage of $\mathcal{B}$, and does not yield a sample in which any networks from $\mathcal{B}$ are over-represented.

Finally, using the stopping rule yielded a random sample of 1,000 bipartite networks in a practical amount of time: about 1 second for the Finches and Women, and about 77 seconds for the Senate, using an Apple M1 Max processor with 64GB of RAM.

\section{Discussion}
Researchers performing statistical analysis of a bipartite network need to draw a random sample of all bipartite networks with the same degree sequences as the observed network. Efficient `trade' algorithms exist for sampling bipartite networks through repeated randomization trades \citep{strona2014fast,godard2022fastball}. Although using these algorithms on contemporary computers can perform trades quite fast, and future algorithms executed using future computing technology will enable trades to be performed even faster, it is still necessary to know how many trades needed to ensure the sample is random. In this paper, I proposed a stopping rule that focuses on the distribution of the sampled networks' distances from a starting network, and stops performing trades once this distribution stabilizes. For a wide range of degree sequences, this stopping rule ensures with high probability (87.5\% on average, 92\% adjusted for type-I error) that the sample is a random sample of all bipartite networks with the given degree sequence. Additionally, in three larger empirical examples, the stopping rule performs as expected, and returns a sample in a practical amount of time.

The stopping rule offers a solution to a practical challenge encountered by researchers: determining when a sample of bipartite networks is random, regardless of how many trades this might require. Importantly, this challenge cannot be solved by greater computational power alone. Greater computational power may allow trades to be performed more quickly, thereby enabling faster sampling and sampling of larger bipartite networks, but researchers still require a rule for deciding \textit{how many} trades to perform.

These analyses do not directly answer a related question -- how many trades are required, or what is the mixing time of a trade algorithm -- but they do offer some insight. In their provided code, \cite{strona2014fast} implicitly suggested that $5T$ trades was sufficient, where $T$ is the number of top nodes, but they did not provide a rationale or evidence for this suggestion. Was their suggestion reasonable?

For the small networks described in section \ref{sec:small}, an average of about 4 trades per top node were required to obtain a random sample ($M = 4.11, SD = 0.50$). These results suggest that a $5T$ rule-of-thumb may be slightly liberal, but not unreasonable, in small networks. In contrast, for the empirical networks described in section \ref{sec:large}, many more trades were required: $9T$ for Darwin's Finches, $7T$ for Southern Women, and $11T$ for the US Senate. These results suggests that in larger networks a $5T$ rule-of-thumb may be too conservative, and that more trades are needed to achieve a random sample. More generally, contrary to the suggestion from \cite{strona2014fast}, these results suggest that there is not a linear relationship between the number of top nodes and the number of trades required. While the results presented above do not provide definitive guidance on how many trades to perform, the stopping rule provides a way to determine how many trades to perform in a specific case, as well as a tool for further research on trade algorithm mixing times.

These findings must be interpreted in light of two computational limitations. First, because $|\mathcal{B}|$ grows very rapidly as a function of network size and in most cases is unknown \citep{barvinok2010number}, it is only possible to evaluate the randomness of a sample for relatively small bipartite networks. Therefore, although the results presented in section \ref{sec:large} are consistent with the resulting samples being random, their randomness cannot be evaluated directly. Second, although the stopping rule provides confidence that the resulting sample is random, this confidence comes at the cost of greater time and memory requirements. In terms of time, using the stopping rule takes longer than simply performing trades because it requires periodically computing $D$ and performing a AD or DTS test. In terms of memory, using the stopping rule requires requires holding all the sampled networks in memory as they are incrementally randomized and compared in terms of their distance from the starting network. However, future research may use the stopping rule to investigate the mixing time of trade algorithms, thereby providing guidance on choosing a number of trades in advance and eliminating the need for the stopping rule itself.

The analysis of bipartite networks and their projections has become common in a range of contexts \citep[e.g.,][]{breiger1974duality,neal2014backbone,gotelli2000null,newman2001scientific}. The statistical analysis of these networks typically requires building null models based on random samples of bipartite networks with the same degree sequence as an observed network. Although many methods exist for generating samples of bipartite networks with fixed degree sequences, it is difficult to be sure that these samples are random, that is, that they uniformly sample from all bipartite networks with the given degree sequences. The proposed stopping rule offers a way to be reasonably sure the sample is random, while simultaneously ensuring that the randomization algorithm runs only as long as necessary to ensure the sample's randomness. It may also offer a way to more precisely estimate the mixing time of trade randomization algorithms. Therefore, it offers a useful tool for statistically analyzing bipartite networks, and for investigating network randomization algorithms.

\section*{Data availability statement}
Supplementary information, and all data and code necessary to replicate these results, are available at \url{https://osf.io/rmwpz/}. 

\section*{Disclosures}
The author declares no conflicts of interest.

\end{document}